\documentclass[11pt,english]{article}
\usepackage[T1]{fontenc}
\usepackage[latin1]{inputenc}

\makeatletter
\newcommand{\lyxaddress}[1]{
\par {\raggedright #1
\vspace{1.4em}
\noindent\par}
}

\title{The Semigroup of a Word}
\author{Peter M. Higgins \& Norman R. Reilly}
\date{}
\def\qed{\quad\vrule height4.17pt width4.17pt depth0pt}

\makeatother

\usepackage{babel}
\begin{document}

\title{\textbf{Orthodox semigroups and permutation matchings}}

\author{Peter M. Higgins }

\maketitle

\lyxaddress{\centerline{Dept of Mathematical Sciences, University of Essex}\centerline{Wivenhoe Park, Colchester, UK, CO4 3SQ}\centerline{email: peteh@essex.ac.uk}}
\begin{abstract}
We determine when an orthodox semigroup $S$ has a permutation that
sends each member of $S$ to one of its inverses and show that if
such a permutation exists, it may be taken to be an involution. In
the case of a finite orthodox semigroup the condition is an effective
one involving Green's relations on the combinatorial images of the
principal factors of $S$. We also characterise some classes of semigroups
via their permutation matchings. 
\end{abstract}

\section{Introduction and Background }

In {[}2{]} the author introduced the topic of \emph{permutation matchings,
}which are permutations on a regular semigroup that map each element
to one of its inverses. The main results of this note are in Section
3 where we characterise the class of orthodox semigroups that have
a permutation matching and show that every orthodox semigroup with
a permutation matching possesses an involution matching. In Section
2 we detail how some classes of semigroups may be characterised by
their permutation matchings.

Following the text of Howie {[}4{]} and the author's {[}3{]} we denote
the set of idempotents of a semigroup $S$ by $E(S)$. We shall write
$(a,b)\in V(S)$ if $a$ and $b$ are mutual inverses in $S$ and
denote this as $b\in V(a)$ so that $V(a)$ is the set of inverses
of $a\in S$. We extend the notation for inverses to sets $A$: $V(A)=\bigcup_{a\in A}V(a)$.
Standard results on Green's relations, particularly those stemming
from Green's Lemma, will be assumed (Chapter 2 of {[}4{]}, specifically
Lemma 2.2.1) and indeed basic facts and definitions concerning semigroups
that are taken for granted in what follows are all to be found in
{[}3,4{]} . Recall that a regular semigroup $S$ is \emph{orthodox}
if $E(S)$ forms a subsemigroup of $S$. We say that a semigroup $S$
is \emph{combinatorial }(or \emph{aperiodic}) if Green's ${\cal H}$-relation
on $S$ is trivial. A completely $0$-simple combinatorial semigroup
is known as a $0$-\emph{rectangular band}. When $S$ is a finite
semigroup, we shall write $a^{\omega}$ for the unique idempotent
power of $a\in S$. We shall also write $a^{\omega+1}$ to denote
$a^{\omega}a$ and let $a^{\omega-1}$ stand for the least positive
power $a^{k}$ of $a$ such that $a^{k+1}=a^{\omega}$. Following
{[}1{]} we say that a member $a\in S$ is \emph{co-regular} if there
exist $b\in S$ such that $a=aba=bab$ and $S$ is co-regular if all
of its members have that property. Co-regularity of $a$ is equivalent
to saying that $a$ is self-inverse (which is the term we shall henceforth
use) for if that is the case then $a$ satisfies $a=a^{3}$ and we
may take $b=a$ in order to satisfy the definition of co-regularity.
On the other hand given that $a$ is co-regular then we have by the
given equations that $ab=abab$ and so
\[
(a^{2}=abab=ab)\Rightarrow(a^{3}=aba=a).
\]

Let $C=\{A_{i}\}_{i\in I}$ be any finite family of finite sets (perhaps
with repetition of sets). A set $\tau\subseteq\bigcup A_{i}$ is a
\emph{transversal }of $C$ if there exists a bijection $\phi:\tau\rightarrow C$
such that $t\in\phi(t)$ for all $t\in\tau$. We assume Hall's Marriage
Lemma in the form that $C$ has a transversal if and only if \emph{Hall's
Condition }is satisfied, which says that for all $1\leq k\leq|I|$,
the union of any $k$ sets from $C$ has at least $k$ members. For
this and related background see, for example, the text {[}6{]}.

\textbf{Definitions 1.1 }Let $S$ be a regular semigroup and let $F=\{f\in T_{S}:f(a)\in V(a)\}$.
We call $F$ the set of \emph{inverse matchings }of $S$. We call
$f\in F$ a \emph{permutation matching }if $f$ is a permutation of
$S$; more particularly $f$ is an \emph{involution matching} if $f^{2}=\iota$,
the identity mapping.

We shall often denote a matching simply by $(\cdot)'$, so that the
image of $a$ is $a'$. When discussing an involution $f$ we may
sometimes write $a\leftrightarrow b$ to indicate that $f(a)=b$ and
$f(b)=a$. In general the inverse $f^{-1}$ of a permutation matching
$f$ is also a permutation matching, and so if a semigroup has a unique
permutation matching $f$, then $f$ must necessarily be an involution
matching. Permutation matchings are not however closed under composition
of permutations. 

When using the following theorem we work with the family of subsets
of $S$ given by $V=\{V(a)\}_{a\in S}$. The members of $V$ may have
repeated elements---for example $S$ is a rectangular band if and
only if $V(a)=S$ for all $a\in S$. However, for the purposes of
the next result we consider the members of $V$ to be marked by the
letter $a$, so that $V(a)$ is an unambiguous member of $V$ (strictly,
we are using the pairs $\{a,V(a)\},$ $(a\in S)$). We summarise some
results of {[}2{]}.

\textbf{Proposition 1.2 }For a finite regular semigroup $S$ the following
are equivalent:

(i) $S$ has a permutation matching;

(ii) $S$ is a transversal of $V=\{V(a)\}_{a\in S}$;

(iii) $|A|\leq|V(A)|$ for all $A\subseteq S$. 

\textbf{Remark 1.3 }In general, a finite regular semigroup may have
no permutation matchings. A minimal example is given by the $7$-element
orthodox $0$-rectangular band $B=\{(i,j):1\leq i\leq2,\,1\leq j\leq3\}\cup\{0\}$,
where $E(B)=\{(1,2),(1,3),(2,1)\}\cup\{0\}$. Then $V\{(2,2),(2,3)\}=\{(1,1)\}$
and so Hall's Condition is violated, and by Proposition 1.2, $B$
has no permutation matching. The class of finite regular semigroups
with a permutation matching is closed under the taking of direct products
but not of homomorphic images or regular subsemigroups {[}2, Proposition
1.5.{]} 

\textbf{Theorem 1.4 }{[}see 2, Theorem 1.6{]}\textbf{ }The following
are equivalent for a finite regular semigroup $S$:

(i) $S$ has a permutation matching;

(ii) $S$ has a permutation matching that preserves the ${\cal H}$-relation
(meaning that $\alpha{\cal H}\beta\Rightarrow\alpha'{\cal H}\beta'$);

(iii) each principal factor $D_{a}\cup\{0\}$ $(a\in S)$ has a permutation
matching;

(iv) each $0$-rectangular band $B=D_{a}\cup\{0\}/{\cal H}$$(a\in S)$
has a permutation matching.

\textbf{Remarks 1.5} We may re-write the proof of Theorem 1.4 as given
in {[}2{]}, replacing the word `permutation' by `involution' to recover
the implications ((i) $\Leftrightarrow$ (iii)) $\Leftarrow$ ((ii)
$\Leftrightarrow$ (iv)). However the missing forward implication
has not been proved and remains an open question.

\emph{Proof} (i) $\Leftrightarrow$ (iii) Suppose that $':S\mapsto S$
is an involution matching of $S$ and let $D$ denote a ${\cal D}$-class
of $S$. For each $a\in D$ we have $V(a)\subseteq D$ so that $'|_{D}$
is an involution of $D$; extending this by $0\mapsto0$ gives an
involution matching for the principal factor $D\cup\{0\}$. Conversely,
if each principal factor has an involution matching then the union
of these matchings over the set of ${\cal D}$-classes yields an involution
matching of $S$.

(ii) $\Rightarrow$(iv) The given involution matching $(\cdot)'$
of $S$, when restricted to a ${\cal D}$-class $D$, defines an involution
matching on the principal factor $D\cup\{0\}$. In addition we are
told that $(\cdot)'$ preserves the ${\cal H}$-relation and so $(\cdot)'$
induces an involution matching on $B$.

(iv) $\Rightarrow$ (ii) Let $D$ be an arbitrary ${\cal D}$-class
of $S$.  The given involution matching of $B$ induces an involution
on the set of ${\cal H}$-classes within $D$ in such a way that $H_{1}\leftrightarrow H_{2}$
then implies that each $a\in H_{1}$ has a unique inverse $a'\in H_{2}$.
The mapping defined by $a\leftrightarrow a'$ then defines an involution
of $H_{1}\cup H_{2}$. Taking the union of these involutions over
all such pairs $(H_{1},H_{2})$ then gives a required involution matching
of $D$. (Note that $H_{1}=H_{2}$ is possible, in which case $H_{1}$
is a group and the involution $a\leftrightarrow a'$ is the unique
involution matching of that group.). 

(ii) $\Rightarrow$(i) as a special case.

\textbf{Theorem 1.6 }{[}3, Theorem 1.4.18{]} For a regular semigroup
$S$, $V(E^{n})=E^{n+1}$. In particular for an orthodox semigroup,
$V(E)=E^{2}=E$. 

\textbf{Theorem 1.7 }{[}4, see Theorems 6.2.4 and 6.2.5{]} Let $S$
be a semigroup. Then $S$ is orthodox if and only if $\{V(a):a\in S\}$
forms a partition of $S$. Let $\gamma$ be the relation on the semigroup
$S$ whereby $a\gamma b$ if and only if $V(a)=V(b)$. The relation
$\gamma$ is an equivalence relation for any semigroup $S$. If $S$
is orthodox then $\gamma$ is the least inverse semigroup congruence
on $S$. 

We close this section with several observations based on these results
that will be invoked in Section 3.

\textbf{Corollary 1.8 }For the orthodox semigroup $S$, the set of
inverse sets $I=\{V(a):a\in S\}$ equals the set of classes of the
least inverse congruence $\gamma$ of $S$. Moreover for any $a\in S$,
$V(V(a))=a\gamma$.

\emph{Proof }For any semigroup $S$ the relation $\gamma$ whereby
$a\gamma b$ if $V(a)=V(b)$ is an equivalence relation. If $(a,b)\in V(S)$
then $a\gamma\subseteq V(b)$. It follows that if $S$ is regular,
every $\gamma$-class is contained in $V(b)$ for some $b\in S$.
If we assume further that $S$ is orthodox we have by Theorem 1.7
that the classes of inverses $I$ also partition $S$. Since each
$\gamma$-class is contained in some class of inverses it follows
that each class of inverses $V(a)$ is the union of $\gamma$ classes.
However if $b,c\in V(a)$ then $a\in V(b)\cap V(c)$ whence $V(b)=V(c)$,
which is to say $b\gamma c$. Hence each class of inverses $V(a)$
consists of exactly one $\gamma$ class. Therefore the two partitions
of $S$ are identical: $\{a\gamma:\,a\in S\}=\{V(a):\,a\in S\}$.

For the second statement take any $b\in V(a)$. Then since $a\in V(b)$
it follows that $a\in V(V(a))$. For any $c\in a\gamma$ we also have
$c\in V(b)$ and so it follows that $a\gamma\subseteq V(V(a))$. Indeed,
since $b\in V(a)$ was arbitrary, it follows that $a\gamma\subseteq\cap_{b\in V(a)}V(b)$.
However since $S$ is orthodox, the distinct classes in this intersection
are pairwise disjoint. Since $a\gamma\neq\emptyset$ it follows that
all members of this intersection are equal and so $V(V(a))$ is itself
a single class of inverses. Since $a\gamma$ is a member of $I$ that
is contained in $V(V(a))$, it follows that $V(V(a))=a\gamma$, as
required. $\qed$

\textbf{Remark} \textbf{1.9}: It follows from Theorem 1.7 and Corollary
1.8 that for any member $a$ of an orthodox semigroup $S$, either
$a\gamma\cap V(a)=\emptyset$ or $a\gamma=V(a)$, the latter occurring
exactly when $a=a^{3}$. In this case all members $b\in V(a)=a\gamma$
are self-inverse for it follows that $b\gamma=a\gamma$ and so $b\in V(b)$
and $b=b^{3}$ also.

\textbf{Corollary 1.10} For an orthodox semigroup $S$ the mapping
$V:V(a)\mapsto V(V(a))$ is an involution on the set $I=\{V(a):a\in S\}$.
Moreover the fixed points of this involution are exactly the classes
$V(a)$ where $a=a^{3}$.

\emph{Proof }By Corollary 1.8, any member of $I$ has the form $a\gamma$
$(a\in S)$ and conversely each class $a\gamma\in I$. Hence $V(V(a))=a\gamma$
shows that $V$ is a mapping from $I$ into $I$. Moreover $V(V(a\gamma))=V(V(a))=a\gamma$
shows that $V$ does indeed define an involution on $I$. By Remark
1.9, $V(a\gamma)=a\gamma$ if and only if $a$ is self-inverse, which
occurs if and only if every member of $a\gamma$ is self-inverse.
Hence $V(a)$ is fixed by our involution if and only if $a=a^{3}$.
$\qed$

\section{Characterisation of classes of regular semigroups by their permutation
matchings}

In this section we show how classes of finite regular semigroups can
be characterised by their permutation matchings. We begin with \emph{completely
regular }semigroups, which are those that are unions of groups.

\textbf{Theorem 2.1 }Let $S$ denote a finite regular semigroup.

(i) $S$ is completely regular if and only if the function $f(x)=x^{\omega-1}$
$(x\in S)$ is a permutation matching;

(ii) $S$ is completely simple if and only if $f(x)=x^{\omega-1}(xyx)^{\omega}$
$(x,y\in S)$ is a permutation matching; 

(iii) $S$ is a group if and only if $f(x)=y^{\omega}x^{\omega-1}y^{\omega}$
$(x,y\in S)$ is a permutation matching;

(iv) $S$ satisfies $x=x^{k+2}$ $(k\geq1)$ if and only if $f(x)=x^{k}$
$(x\in S)$ is a permutation matching; in particular $S$ is co-regular
if and only if the identity function is a permutation matching;

(v) $S$ is a rectangular band if and only if every permutation of
$S$ is a permutation matching.

\textbf{Remark}: In (ii) and (iii) the function $f$ is in general
a function of the two variables $x$ and $y$ but, under the hypothesis
of each part, the output of $f$ depends on $x$ alone. We record
proofs for just the first two statements, the others being straightforward
and similar.

\emph{Proof }(i) Given that $S$ is a union of groups, for any $x\in S$
we have $xx^{\omega-1}x=x^{\omega}x=x$ as $x^{\omega}$ is the identity
of the group $H_{x}$; similarly we obtain $x^{\omega-1}xx^{\omega-1}=x^{\omega-1}x^{\omega}=x^{\omega-1}$.
Since $x$ and $y=x^{\omega-1}$ each have no other inverse in the
group $H_{x}=H_{y}$ it follows that $f(x)$ is a permutation matching,
indeed $f(x)$ is an involution matching. Conversely, given that $f(x)$
is a permutation matching of $S$ we have that $(x,x^{\omega-1})$
are mutual inverses in $H_{x}$, which is therefore a group, and so
$S$ is a union of groups.

(ii) Given that $S$ is completely simple it follows that for any
$x,y\in S$ we have that $(xyx)^{\omega}=x^{\omega}$, the identity
element of the group $H_{x}$. Hence $f(x)=x^{\omega-1}(xyx)^{\omega}=x^{\omega-1}$
and as in (i) we have that $x$ and $x^{\omega-1}$ are mutually inverse
in the group $H_{x}$ and $f(x)$ is an involution matching. Conversely
suppose that $f(x)$ is a permutation matching. If $S$ had more than
one ${\cal D}$-class then there would exist ${\cal D}$-classes,
$D_{1}$ and $D_{2}$ such that $D_{1}<_{{\cal {\cal D}}}D_{2}$.
However if we then take $x\in D_{2}$ and $y\in D_{1}$ we get that
$f(x)\not\in D_{2}$, contradicting our assumption that $f(x)\in V(x)$.
Hence $S$ is completely simple.$\qed$ 

\textbf{Theorem 2.2 }An inverse semigroup $S$ has a unique permutation
matching, which is an involution. Conversely, any orthodox semigroup
$S$ with a unique permutation matching is an inverse semigroup. However,
there exists a $5$-element combinatorial semigroup that is not an
inverse semigroup, which has a unique permutation matching.

\emph{Proof} For an inverse semigroup $S$, clearly the mapping whereby
$a\mapsto a^{-1}$ is the unique permutation matching of $S$ and
is an involution. Conversely, suppose that $S$ is orthodox and has
a unique permutation matching $f$. If $S$ is not inverse then there
exists two distinct idempotents $e,f\in E=E(S)$ such that $e{\cal G}f$,
where ${\cal G}$ denotes either of Green's relations ${\cal L}$
or ${\cal R}$. Since $S$ is orthodox, it follows by Theorem 1.6
that $f(E)=E$ and $f(S\setminus E)=S\setminus E$. However these
equations now allow us to construct two distinct permutation matchings
$g,h$ of $S$ as follows. Put $g|_{(S\setminus E)}=h_{|(S\setminus E)}=f_{|(S\setminus E)}$
and put $g_{|E}=\iota_{|E}$ the identity mapping on $E$, while $h$
is identical to $g$ except that $h(e)=f$ and $h(f)=e$. This contradicts
the uniqueness of matchings and so it follows that if $S$ is orthodox
with a unique permutation matching then $S$ is inverse and that matching
is the standard involution by inverses. 

Next consider the $5$-element combinatorial semigroup $S=D\cup\{0\}$
where the $2\times2$ ${\cal D}$-class $D=\{(1,1),(1,2),(2,1),(2,2)\}$
consists of elements that are all idempotent except $(1,2)$. This
semigroup is regular but not orthodox as the idempotent $(2,1)$ is
the unique inverse of the non-idempotent $(1,2).$ However $S$ has
a unique permutation matching \textbf{$f$}: for any permutation matching
$f$ we necessarily have $(1,2)\mapsto(2,1)$; we can now complete
an involution matching $f$ by saying that under $f$, $(2,1)\mapsto(1,2)$
while the three idempotents $0,\,(1,1)$, and $(2,2)$ are fixed by
$f$. What is more, by inspection we see there is no alternative to
this definition for $f$. Therefore $S$ possesses a unique permutation
matching but $S$ is not an inverse semigroup. $\qed$ 

\textbf{Remark 2.3 }We observe that if $S$ is any semigroup that
is not inverse but has a unique permutation matching $f$, then $S$
contains the ${\cal D}$-class structure of $D$ as given in the example
of the previous proof. To see this, note that by the argument of the
proof of Theorem 2.2, the uniqueness of $f$ ensures that $f(E)\not\neq E$
so there exists $e\in E=E(S)$ and $a\in S\setminus E$ such that
$(e,a)\in V$. Hence $(e,a)\not\in{\cal G}$ in $S$ and so $ea=f\in E$
with $e{\cal R}f$ and $ae=g\in E$ with $e{\cal L}g$ and $e,f,g$
are pairwise distinct. Then $gf=ae\cdot ea=aea=a$ and so $\{e,f,g,a\}$
has the ${\cal D}$-class structure of $D$ above, meaning that $e,f,g\in E$,
$a\not\in E$ and $a{\cal R}g{\cal L}e{\cal R}f{\cal L}a$ and $gf=a,\,fg\not\in D$,
$(a,e)\in V(S)$.

\section{Orthodox semigroups}

\textbf{Theorem 3.1} An orthodox semigroup $S$ has a permutation
matching if and only if $|V(V(a)|=|V(a)|$ for all $a\in S$. If orthodox
$S$ has a permutation matching then in fact $S$ possesses an involution
matching.

\emph{Proof} We prove that the condition that $S$ has a permutation
matching implies the stated condition on cardinalities of sets of
inverses and that that condition in turn allows the construction of
an involution matching for $S$. From this follows the claims in the
statement of the theorem. 

To this end suppose that $S$ has a permutation matching $(\cdot)'$.
By Corollaries 1.8 and 1.10 the mapping whereby $V(a)\mapsto V(V(a))$
is an involution on the set $I=\{V(a):a\in S\}=\{a\gamma:a\in S\}$.
For any $a\gamma\in I$ we have $(a\gamma)'\subseteq V(a)$ and since
$(\cdot)'$ is one-to-one, it follows that $|a\gamma|\leq|V(a)|$.
Similarly $(V(a))'\subseteq V(V(a))=a\gamma$ by Corollary 1.10, and
so $|V(a)|\leq|a\gamma|$. It follows that $|a\gamma|=|V(a)|$. Since
any pair of mutual images under $V:I\rightarrow I$ has the form $(a\gamma,V(a\gamma)=V(a))$
it follows that the members of each pair of mutual images under the
involution $V:I\rightarrow I$ are equi-cardinal, which is to say
that $|V(V(a)|=|V(a)|$ for all $a\in S$. 

Let the classes of $I=\{V(a):a\in S\}$ that are fixed by the involution
$V$ of $I$ be $\{I_{i}\}_{i\in R}$ for some index set $R$, and
write the remaining classes of $I$ in ordered pairs $\{(J_{i},K_{i})\}_{i\in T},$
for some index set $T$, where $J_{i}\leftrightarrow K_{i}$ under
the involution $V$. Recall that $J_{i}\cap K_{i}=\emptyset$ and
each pair $(a,b)\in J_{i}\times K_{i}$ is a pair of mutual inverses.
The stated condition on cardinalities of inverse sets says that $|J_{i}|=|K_{i}|$
for each $i\in T$. We may now define an involution matching $f$
on $S$ as follows. Since the members of each pair of elements of
any $I_{i}$ are mutually inverse, we may take $f|_{I_{i}}$ to be
any involution mapping of $I_{i}$ (for example, the identity mapping
on $I_{i}$). Since $J_{i}$ and $K_{i}$ are equi-cardinal, for each
$J_{i}$ we may let $f|_{J_{i}}$ be any bijection of $J_{i}$ onto
$K_{i}$ and let the $f|_{K_{i}}$ be the corresponding inverse bijection.
The mapping $f:S\rightarrow S$ is then the union of all these restriction
mappings on the members of $I$. By construction $f$ is an involution
that maps each member of $S$ to one of its inverses and so $f$ is
an involution matching of $S$, as required. $\qed$ 

We next use Theorem 3.1 to determine when a finite orthodox semigroup
$S$ has a permutation matching and to effectively find all such matchings
when they exist. The result hinges on the special case where $S$
is a finite orthodox $0$-rectangular band so to this end let $S$
denote a finite $m\times n$ orthodox $0$-rectangular band (meaning
the number of ${\cal R}$- and ${\cal L}$-classes in $D=S\setminus\{0\}$
is $m$ and $n$ respectively) with band of idempotents $E(S)=B$.
The non-zero members of $S$ can therefore be taken to be the set
of ordered pairs $S\setminus\{0\}=\{(i,j):\,1\leq i\leq m,\,1\leq j\leq n\}$.
This structure for $S$ is assumed in Lemmas 3.2, 3.3 and 3.4.

\textbf{Lemma 3.2 }The set of inverses $V(e)$ of $e\in B$ is a maximal
rectangular subband of $S$. 

\emph{Proof }Consider the minimum inverse congruence $\gamma$ on
$S$. Then $(e\gamma)^{2}=e^{2}\gamma=e\gamma$ and so from Remark
1.9 we obtain $(V(e))^{2}=V(e)$, which is to say that $V(e)$ is
a subsemigroup of $S$ and by Theorem 1.6 it follows that $V(e)$
is a band. For any $a$ such that $(a,e)\in V$ we have that $e\in V(e)\cap V(a)$
and so by Theorem 1.7 we infer that $V(a)=V(e)$. It follows that
for any $a,b\in V(e)$ we have $b\in V(a)$ also and so $aba=a$.
Therefore $V(e)$ is a rectangular subband of $S$. 

Next let $V(e)\subseteq U$, where $U$ is a rectangular subband of
$S$. Let $e{\cal G}f$ in $U$, whence $e{\cal G}f$ in $S$. Then
$f\in V(e)\cap V(f)$ so that $V(e)=V(f)$. Take any $u\in U$. Since
$U$ consists of a single ${\cal D}$-class and $V(e)\subseteq U$
it follows that there exists $f\in U$ such that, in $U$, $e{\cal R}f{\cal L}u$.
By the previous argument this gives first that $f\in V(e)$ and then
in turn $u\in V(e)$. Therefore $U\subseteq V(e)$ and so we conclude
that $U=V(e)$, which is to say that $V(e)$ is indeed a maximal rectangular
subband of $S$. $\qed$

Denote the pairwise distinct maximal rectangular subbands formed by
the sets of inverses $V(e)$ $(e\in B\setminus\{0\})$ by $V(e_{1}),V(e_{2}),\cdots,V(e_{k})$
say for suitable fixed representatives $e_{i}$ of each of the rectangular
subbands $V(e).$ 

\textbf{Lemma 3.3 }Each ${\cal G}$-class of $D=S\setminus\{0\}$
meets exactly one of the maximal rectangular subbands $V(e)$.

\emph{Proof }By symmetry it is enough to take the case ${\cal G}={\cal R}$
so let $a\in D$, consider the class $R_{a}$ and take any $e\in B\cap R_{a}$
so that $R_{a}$ meets the maximal rectangular subband $V(e)$. Suppose
that $V(e_{i})$ is a maximal rectangular subband of $D$ such that
$R_{a}\cap V(e_{i})\neq\emptyset$. Then any $f\in R_{a}\cap V(e_{i})$
is idempotent, $e{\cal R}f$ and so $f\in V(e)$. However we then
have $f\in V(e_{i})\cap V(e)$ so that $V(e_{i})=V(e)$, thus completing
the proof. $\qed$

We may now arrange the `egg-box' diagram of $D$ so that the maximal
rectangular subbands, $U_{1}=V(e_{1}),U_{2}=V(e_{2}),\cdots$ form
a leading diagonal of the diagram as follows. We may list the ${\cal {\cal L}}$-classes
of $D$ (and similarly the ${\cal R}$-classes) by first listing all
the ${\cal {\cal L}}$-classes of $D$ that meet $U_{1}$, then of
$U_{2}$ and so on through to those of $U_{k}$, as, by Lemma 3.3,
these orderings of the ${\cal L}$- and the ${\cal R}$-classes are
well-defined. Let $U$ denote the set of members of $D$ consisting
of all the intersections of the ${\cal L}$-classes in $S$ that meet
$U_{i}$ with the ${\cal R}$-classes in $S$ that meet $U_{i}$.
Clearly $U_{i}\subseteq U$ but the reverse inclusion is also true
for take any $u\in U$ so that $\{u\}=R_{f}\cap L_{g}$ say where
$f,g\in U_{i}=V(e_{i})$. Since $V(e_{i})$ is a subsemigroup of $S$
it follows that $fg\in U_{i}$; in particular $fg\neq0$ so that $u=fg\in U_{i}$.
Therefore the union of intersections of the ${\cal L}$- and ${\cal R}$-classes
of $U_{i}$ forms the maximal rectangular subband $U_{i}$ itself.

Let the number of ${\cal L}$- and ${\cal R}$-classes of $V(e_{i})$
be denoted by $m_{i}$ and $n_{i}$ respectively $(1\leq i\leq k)$.
Define a mapping 
\[
\phi:S\setminus\{0\}\rightarrow\{(i,j):\,1\leq i,j\leq k\}
\]
where $\phi(a)=(i,j)$ if $V(e_{i})$ and $V(e_{j})$ are the maximal
rectangular subbands of $S$ that meet $R_{a}$ and $L_{a}$ respectively. 

\textbf{Lemma 3.4 }Let $a\in S\setminus\{0\}$ with $\phi(a)=(i,j)$.
Then

(i) $V(a)=\{b\in S:\,\phi(b)=(j,i)\}=V(\{a\in S:\,\phi(a)=(i,j)\})$;

(ii) Let the set defined in (i) be denoted by $V_{i,j}$. Then $V(V_{i,j})=V_{j,i}$.

(iii) $|V(a)|=m_{j}n_{i}$.

\emph{Proof} (i) and (ii). In general, $b\in V(a)$ if and only if
there exist $e,f\in E$ such that $a{\cal R}e{\cal L}b{\cal R}f{\cal {\cal L}}a$.
Since $\phi(a)=(i,j)$ this implies that $e\in V(e_{i}),\,f\in V(e_{j})$,
whence $\phi(b)=(j,i)$. Conversely if $\phi(b)=(j,i)$ then for some
$e\in V(e_{j})$ and $f\in V(e_{i})$ we have $f{\cal L}b{\cal R}e$.
Let $\{g\}=L_{f}\cap R_{a}$ and $\{h\}=R_{e}\cap L_{a}$. Now $a{\cal R}k$
for some $k\in V(e_{i})$ so that $g=kf\in V(e_{i})$ and similarly
$h\in V(e_{j})$. Therefore we have $a{\cal R}g{\cal L}b{\cal R}h{\cal L}a$
and since $g,h\in E(S)$ we conclude that $b\in V(a)$. It now follows
that for any $a_{1},a_{2}\in D$, $V(a_{1})=V(a_{2})$ if and only
if $\phi(a_{1})=\phi(a_{2})$, whence the second equality in (i) now
follows, while (ii) is a re-statement of the second equality in (i).

(iii) By (i) we have $|V(a)|=|\{b\in S:\,\phi(b)=(j,i)\}|$. To say
that $\phi(b)=(j,i)$ means that $R_{b}\cap V(e_{j})\neq\emptyset$
and $L_{b}\cap V(e_{i})\neq\emptyset$. This in turn gives $m_{j}$
choices for $R_{b}$ and $n_{i}$ choices for $L_{b}$. Since $b$
is determined by such a pair of choices, which can be made independently
of each other, it follows that $|V(a)|=m_{j}n_{i}$, as claimed. $\qed$

\textbf{Definition 3.5 }Let $U_{1}$ and $U_{2}$ be finite rectangular
bands, let $m_{i}$ (resp. $n_{i})$ be the respective number of ${\cal R}$-classes
and ${\cal L}$-classes of $U_{i}$ $(i=1,2)$. We shall say that
$U_{1}$ and $U_{2}$ are \emph{similar }if $\frac{m_{1}}{n_{1}}=\frac{m_{2}}{n_{2}}$. 

\textbf{Theorem 3.6 }Let $S$ be a finite orthodox $0$-rectangular
band. Then $S$ has a permutation matching if and only if the maximal
rectangular subbands of $D=S\setminus\{0\}$ are pairwise similar. 

\emph{Proof} Using the notation of Lemma 3.4, we note that the sets
$V_{i,j}$ partition $D$ into blocks and, by Lemma 3.4(i) and (ii),
two members $a,b\in D$ lie in the same block $V_{i,j}$ if and only
if $V(a)=V(b)$. Suppose that $S$ has a permutation matching $f$.
It follows from Lemma 3.4(ii) that $f_{|V_{i,j}}$ is a bijection
onto $V_{j,i}$ and in particular that $|V_{i,j}|=|V_{j,i}|$ for
any pair $(i,j)$ $(1\leq i,j\leq k)$. Hence from Lemma 3.4(iii)
we now infer: 
\begin{equation}
(|V_{i,j}|=|V_{j,i}|)\Leftrightarrow(m_{j}n_{i}=m_{i}n_{j})\Leftrightarrow\big(\frac{m_{i}}{n_{i}}=\frac{m_{j}}{n_{j}}\big)
\end{equation}
which is equivalent to the condition that the members of each pair
of maximal rectangular subbands $U_{i}$ and $U_{j}$ of $D$ are
similar. 

Conversely, given that the members of each such pair of rectangular
subbands in $D$ are similar, it follows that (1) holds for each pair
$(i,j)$. It now follows from Lemma 3.4 and Theorem 3.1 that $S$
has a permutation matching. Indeed all such permutation matchings
may be constructed as the union of any chosen bijections $':V_{i,j}\rightarrow V_{j,i}$
$(1\leq i,j\leq k)$. The involution matchings correspond to the case
where the members of each pair of bijections between sets $V_{i,j}$
and $V_{j,i}$ are mutual inverses. $\qed$

\textbf{Theorem 3.7 }Let $S$ be a finite orthodox semigroup. Then
$S$ has a permutation matching if and only if for each $0$-rectangular
band $D_{a}\cup\{0\}/{\cal H}$$(a\in S)$ the maximal rectangular
subbands are pairwise similar.

\emph{Proof }This follows from Theorem 1.4 (i) $\Leftrightarrow$(iv)
and Theorem 3.6.

It remains an open question as to whether there exists a finite regular
semigroup (necessarily non-orthodox) that has a permutation matching
but no involution matching. In particular, it was shown in {[}2, Theorem
2.12{]} that any finite full transformation semigroup possesses a
permutation matching but, as yet, no involution matching has been
identified for $T_{n}$. In {[}5{]} it was shown that $T_{n}$ is
covered by its inverse subsemigroups (also see {[}3, Chapter 6.2{]})
but that fact does not in itself immediately yield an involution matching
for $T_{n}$.

\end{document}